\documentclass[a4paper]{article}

\usepackage{amsmath}
\usepackage{amssymb}
\usepackage{amsthm}
\usepackage{indentfirst}

\makeatletter

\newcommand{\Rmnum}[1]{\uppercase\expandafter{\romannumeral #1}}

\newtheorem{thm}{Theorem}[section]

\newtheorem{lem}[thm]{Lemma}

\theoremstyle{definition}
\newtheorem{rem}[thm]{Remark}

\begin{document}

\title{\textbf{On Estimates of the Number of Collisions for
Billiards in Polyhedral Angles}}
\author{Lizhou Chen\footnotemark}
\date{}
\maketitle

\renewcommand{\thefootnote}{*}
\footnotetext{Institute of Mathematics, Fudan University, Shanghai
200433, China.
\\\hspace*{15pt}E-mail: 031018007@fudan.edu.cn
\\\hspace*{15pt}This work was supported by the Special Funds for Chinese Major
State Basic Research Projects "Nonlinear Science".
}

\noindent{\textbf{Abstract}} {\small{\ We obtain an upper bound of
the number of collisions of any billiard trajectory in a polyhedral
angle in terms of the minimal eigenvalue of a positive definite
matrix which characterizes the angle. Elements of the matrix are
scalar products between the unit normal vectors of faces of the
angle.}


}

\section{Introduction}
Let $H_1,\ldots,H_n$ be $n$ hyperplanes in $m$-dimensional
Euclidean space $\mathbb{R}^m$ passing through the origin and take
a unit normal $\alpha_i$ for each $H_i$, $i=1,\ldots,n$. Assume
the hyperplanes are in general position, i.e.\
$\bigcap_{i=1}^nH_i$ is an ($m-n$)-dimensional plane, or
equivalently, $\alpha_1,\ldots,\alpha_n$ are linear independent.
The polyhedral cone corresponding to $\alpha_1,\ldots,\alpha_n$ is
$$Q=\left\{y\in\mathbb{R}^m\bigm|(y,\alpha_i)\geqslant0,\forall\,i\right\}.$$
Our fundamental object of study is a billiard in $Q$, that is, a
point particle moves with a uniform motion in the interior of $Q$
and has specular (optical) reflections at walls $B_i=H_i\bigcap Q$,
$i=1,\ldots,n$. If a billiard trajectory reaches one of the corners
$B_i\bigcap B_j$, $i\neq j$, its further motion is not defined.

The model is interesting partially because any hard ball system on a
line is isomorphic to a billiard in an appropriate polyhedral
cone---elastic collisions between hard balls correspond to specular
reflections made by the billiard at walls. More generally, all hard
ball systems are isomorphic to semi-dispersing billiards. Estimates
on the number of collisions of billiard trajectories have been
studied for a long time, cf.\ \cite{BFK00}. We present here in a
nutshell only a few selected results rather than a comprehensive
review.

In 1978, Sinai \cite{Sin78} proved the existence of uniform
estimates of the number of collisions of billiard trajectories in a
polyhedral angle. At the same time, he also pointed out that the
smooth version of the result should also hold, that is, when the
polyhedral angle is replaced by a smooth hypersurface with
nonnegative second fundamental form, uniform estimates still exist
in a neighborhood of a point with a condition of linear
independence. Using the same method of Sinai, in 1993, Sevryuk
\cite{Sev94} gave a uniform estimate for billiards in a polyhedral
angle in terms of a geometrical characteristic of the angle. A
milestone is established in 1998 by Burago, Ferleger and Kononenko
\cite{BFK98}, see also \cite{BFK02}, uniform estimates were obtained
for semi-dispersing billiards on arbitrary Riemannian manifolds with
boundaries satisfying a nondegenerate condition.

Now we describe the estimates mentioned above in more details.
Sevryuk introduced the concepts of charge and capacity of a
polyhedral angle in \cite{Sev94}. The charge of the polyhedral
cone $Q$ is defined to be
$$S(Q)=\max_{a\subset Q}\min_{1\leqslant i\leqslant n}\measuredangle(a,H_i),$$
where $0\leqslant\measuredangle(a,H_i)\leqslant{\pi\over2}$
denotes the angle between the ray $a$ emanating from the origin
and the hyperplane $H_i$ and the maximum is taken over all rays
$a$ that pass within $Q$. The charge $\varphi$ of the set of
hyperplanes $H_1,\ldots,H_n$ is the minimum of the charges of
polyhedral cones determined by the hyperplanes. Then
$0<S(Q)\leqslant{\pi\over2}$, $0<\varphi\leqslant{\pi\over2}$. The
capacity $\psi$ of the set of hyperplanes $H_1,\ldots,H_n$ is
$$
0<\min_{a\subset\mathbb{R}^m}\max_{1\leqslant i\leqslant n}
\measuredangle(a,H_i)<{\pi\over2},
$$
if the hyperplanes do not all pass through a common line, and
$\pi\over2$ otherwise. Sevryuk proved $\varphi\leqslant\psi$ and
that arbitrary billiard trajectory in $Q$ has no more than
${\sin^2\varphi\over2}\left({4\over\sin^2\varphi}\right)^{2^{n-1}}-1$
collisions. The estimate of Burago, Ferleger and Kononenko involves
a nondegenerate constant. Following \cite{BFK98}, the nondegenerate
constant of the polyhedral cone $Q$ is
\begin{equation}\label{nondegenerate constant 1}
C=\min_{y\in Q\setminus\left(\bigcap_{j=1}^nB_j\right)}
\max_{1\leqslant i\leqslant n}
\frac{\mbox{dist}(y,B_i)}{\mbox{dist}\left(y,\,\bigcap_{j=1}^nB_j\right)}.
\end{equation}
Then $0<C\leqslant1$. It is proved that any billiard trajectory in
$Q$ has no more than $8\left({1\over C}+2\right)^{2(n-1)}$
collisions.

Consider the positive definite matrix
$\big((\alpha_i,\alpha_j)\big)_{n\times n}$. Let $\lambda_{min}$ be
its minimal eigenvalue. In this paper, we will prove that the number
of collisions of a billiard trajectory in $Q$ does not exceed
$n!\left({4\over\lambda_{min}}\right)^{n-1}$. The case $n=1$ is
obvious and the arguments proceed by induction on $n$---the number
of walls---as in every proof, as far as we know, of finiteness of
the number of collisions for polyhedral billiards when $n>2$.

\section{Estimate by $\lambda_{min}$}
In what follows, we assume the normals $\alpha_1,\ldots,\alpha_n$
span the whole ambient space $\mathbb{R}^m$ (so $m=n$), i.e.\
$\bigcap_{i=1}^nH_i=\{0\}$. Otherwise, we can project the dynamics
to the orthocomplementation of $\bigcap_{i=1}^nH_i$.

Before establishing our estimate for the general case, we would like
to discuss the interesting case $n=2$.

It is well known that unfolding a billiard trajectory inside a wedge
to a straight line yields a sharp bound
$\lceil{\pi\over\theta}\rceil$ for the number of collisions, where
$\theta=\arccos\big(-(\alpha_1,\alpha_2)\big)$ is the angle of the
wedge and $\lceil x\rceil$ is the ceiling function, the smallest
integer not less than $x$, cf.\ \cite{Tab95}. The argument also
shows that if the point particle does not hit the corner of the
wedge at first collision, it will never hit the corner in the
future.

One can take another way as follows in which only the velocity,
rather than the position, of the point particle is concerned.
Suppose the particle moves with unit speed and has suffered $N$
collisions. Let $\mathbf{v}_0,\mathbf{v}_1,\ldots,\mathbf{v}_N$ on
the unit circle be the sequence of velocities. For any $k$,
$1\leqslant k\leqslant N-1$, $\mathbf{v}_k-\mathbf{v}_{k-1}$ is
parallel to some $\alpha_i$, say $\alpha_1$, so
$\mathbf{v}_{k+1}-\mathbf{v}_k$ is parallel to $\alpha_2$. That is
$$\mathbf{v}_k-\mathbf{v}_{k-1}=\|\mathbf{v}_k-\mathbf{v}_{k-1}\|\,\alpha_1,\qquad
\mathbf{v}_{k+1}-\mathbf{v}_k=||\mathbf{v}_{k+1}-\mathbf{v}_k||\,\alpha_2.$$
It is easy to see that the angle $\angle
\mathbf{v}_{k-1}\mathbf{v}_k\mathbf{v}_{k+1}$ is equal to
$\theta$. By elementary geometry, the arc length of
$\widehat{\mathbf{v}_{k-1}\mathbf{v}_{k+1}}$ is $2\theta$ and
hence $2\theta(N-1)<2\pi$. Thus $N<{\pi\over\theta}+1$, equivalent
to $N\leqslant\lceil{\pi\over\theta}\rceil$.

In the general case, we do not intend to find the best estimate for
the number of collisions especially by the method of induction. For
some special cases, sharp bounds may be found as we have seen for
$n=2$. For another example, if $(\alpha_i,\alpha_j)=0$ for
$|i-j|>1$, and $(\alpha_i,\alpha_{i+1})\geqslant-{1\over2}$ for
$i=1,\cdots,n-1$, then the maximal possible number of collisions of
a billiard trajectory is $n(n+1)\over2$, see \cite{Chen}.

Now we proceed to establish our estimate for the general case.

It is convenient to assume the particle moves with unit speed.
Suppose a part of its trajectory has undergone $N$ reflections at
the walls.  Let $\mathbf{v}_0,\mathbf{v}_1,\ldots,\mathbf{v}_N$ on
the unit sphere be the sequence of velocities.

Since $\alpha_1,\ldots,\alpha_n$ are linear independent, they,
perceived as $n$ points in $\mathbb{R}^n$, determine a hyperplane
not passing through the origin. Let $d>0$ be the distance from the
origin to the hyperplane and $\mathbf{e}$ the unit outer normal of
the hyperplane. They are characterized by the equations
\begin{equation}\label{d and e 1}
(\mathbf{e},\alpha_i)=d>0,\,\forall\,i.
\end{equation}
It means that $d$ is the radius of the inscribed ball of the
polyhedral cone $Q$ with the center $\mathbf{e}$ on the unit sphere.
The construction proves the existence of the inscribed ball. Define
the matrix $A=(\alpha_1,\ldots,\alpha_n)$, where $\alpha_i$ are
perceived as column vectors in $\mathbb{R}^n$. Then equation (\ref{d
and e 1}) read
\begin{equation}\label{d and e 2}
\mathbf{e}^T\!A=d(1,\ldots,1).
\end{equation}
Thus
$$
{1\over d}=\left\|(1,\ldots,1)A^{-1}\right\|,\qquad
\mathbf{e}^T\!=d(1,\ldots,1)A^{-1}.
$$
\begin{lem}\label{upper bound of the length L}
Let $L$ be the length of the zigzag line determined by the points
$\mathbf{v}_0,\mathbf{v}_1,\ldots,\mathbf{v}_N$. Then
$$L=\sum_{k=0}^{N-1}\|\mathbf{v}_{k+1}-\mathbf{v}_k\|\leqslant{2\over
d}.$$
\end{lem}
\begin{proof}
From the law of reflection we have
$$
\mathbf{v}_{k+1}-\mathbf{v}_k=\|\mathbf{v}_{k+1}-\mathbf{v}_k\|\,\alpha_{i_k},\quad
k=0,1,\ldots,N-1.
$$
Combining with (\ref{d and e 1}) yields
$$
(\mathbf{v}_{k+1}-\mathbf{v}_k,\,\mathbf{e})=d\,\|\mathbf{v}_{k+1}-\mathbf{v}_k\|,\quad
k=0,1,\ldots,N-1.
$$
Taking the sum over $k$, we obtain
$$L={1\over d}(\mathbf{v}_N-\mathbf{v}_0,\,\mathbf{e})\leqslant{1\over d}\,
\|\mathbf{v}_N-\mathbf{v}_0\|\cdot\|\mathbf{e}\|\leqslant{2\over
d}.$$
\end{proof}
Besides $d$, another constant $\delta$ is involved in our proof. It
has already appeared in the Sinai's original proof \cite{Sin78}.
When $\bigcap_{i=1}^nH_i=\{0\}$, formula (\ref{nondegenerate
constant 1}) reads
\begin{equation}\label{nondegenerate constant 2}
C=\min_{y\neq0\atop y\in Q}\max_{1\leqslant i\leqslant n}
\frac{\mbox{dist}(y,B_i)}{\mbox{dist}(y,0)}=\min_{\|y\|=1\atop
y\in Q}\max_{1\leqslant i\leqslant n}\mbox{dist}(y,B_i).
\end{equation}
We present the definition of $\delta$ similar to formula
(\ref{nondegenerate constant 2}):
$$\delta=\min_{\|y\|=1}\max_{1\leqslant i\leqslant n}\mbox{dist}(y,H_i)=\sin\psi.$$
First of all, $\max_{1\leqslant i\leqslant n}
\mbox{dist}(\,\cdot,H_i)$ is a continuous function. It is positive
everywhere on the unit sphere since $\bigcap_{i=1}^nH_i=\{0\}$. By
compactness of the unit sphere, $\delta>0$.
\begin{thm}
The number of reflections of any billiard trajectory in $Q$ does not
exceed $n!\left({4\over\lambda_{min}}\right)^{n-1}$.
\end{thm}
\begin{proof}
Induction on $n$. The case $n=1$ is trivial and suppose we have
proved the theorem from 1 to $n-1$. Note that in the inductive
hypothesis, the number of walls needs not to be the dimension $m$ of
the configuration space. Now we proceed to prove the theorem for
$n$. At this stage, we may assume $n=m$ as claimed at the beginning
of this section. So we have Lemma \ref{upper bound of the length L}.

Set $N'=(n-1)!\left({4\over\lambda_{min}}\right)^{n-2}$. We need
the fact from linear algebra that the minimal eigenvalue of any
principal submatrix of the positive definite matrix
$\big((\alpha_i,\alpha_j)\big)_{n\times n}=A^T\!A$ is not less
than $\lambda_{min}$. It easily follows from the minimax principle
for eigenvalues, particularly for $\lambda_{min}$:
$$
\lambda_{min}=\min_{\|x\|=1}\|Ax\|^2=\min_{x_1^2+\cdots+x_n^2=1}
\|x_1\alpha_1+\cdots+x_n\alpha_n\|^2.
$$
And by the inductive hypothesis, if a sequence of consecutive
reflections does not involve all the hyperplanes, then the length of
this sequence does not exceed $\left[N'\right]$.

Suppose $\left[N'\right]+1\leqslant N$, thus the first
$\left[N'\right]+1$ reflections involve all the hyperplanes. Hence,
for any $i$, the points
$\mathbf{v_0},\mathbf{v_1},\ldots,\mathbf{v_{\left[N'\right]+1}}$ do
not lie on the same side of the hyperplane $H_i$. Say,
$\mathbf{v_0}$ and $\mathbf{v_{k_i}}$ do not lie on the same side of
$H_i$, then
$\|\mathbf{v}_{k_i}-\mathbf{v}_0\|>\mbox{dist}(\mathbf{v_0},H_i)$.
So the length of the $\left[N'\right]+1$ segments
$$
\sum_{k=0}^{\left[N'\right]}\|\mathbf{v}_{k+1}-\mathbf{v}_k\|>
\max_{1\leqslant i\leqslant n}\mbox{dist}(\mathbf{v_0},H_i)
\geqslant\delta.
$$
It shows that the length of any consecutive $\left[N'\right]+1$
segments of the zigzag line determined by the points
$\mathbf{v}_0,\mathbf{v}_1,\ldots,\mathbf{v}_N$ is bigger than
$\delta$. But Lemma \ref{upper bound of the length L} says that the
length of the whole zigzag line does not exceed $2\over d$. One
obtains
$$N<{2\over{d\delta}}\left(\left[N'\right]+1\right)\leqslant{4\over{d\delta}}N'.$$

It remains to show
${1\over{d\delta}}\leqslant{n\over\lambda_{min}}.$ In fact, $1\over
d$ and $1\over\delta$ are both not bigger than
$\sqrt{n\over\lambda_{min}}$. To this end, we shall use the minimax
principle for minimal eigenvalues and the fact that the minimal
eigenvalue of $AA^T\!$ is the same of $A^T\!A$, which can be seen
from the identity
\begin{equation*}
\begin{split}
\det\left(\lambda I-AA^T\right)&=\det A\ \det\left(\lambda A^{-1}-A^T\right)\\
&=\det\left(\lambda A^{-1}-A^T\right)\ \det A\\
&=\det\left(\lambda I-A^T\!A\right).
\end{split}
\end{equation*}

Taking square of the norm of the two sides of equation (\ref{d and e
2}), one obtains
$$d^2n=\left\|e^T\!A\right\|^2\geqslant\lambda_{min}.$$
On the other hand, $\delta=\max_{1\leqslant i\leqslant
n}\mbox{dist}(y_0,H_i)=\max_{1\leqslant i\leqslant
n}|(y_0,\alpha_i)|$ for some $y_0$ on the unit sphere. Thus
$$\delta^2=\max_{1\leqslant i\leqslant n}|(y_0,\alpha_i)|^2\geqslant
{1\over n}\sum_{i=1}^n|(y_0,\alpha_i)|^2={1\over n}
\left\|y_0^T\!A\right\|^2\geqslant{\lambda_{min}\over n}.$$
\end{proof}
\begin{rem}
The proof also gives another upper bound $\left({4\over d\delta}
\right)^{n-1}$ by setting $N'=\left({4\over d\delta}\right)^{n-2}
$.
\end{rem}
\begin{rem}
Make an observation. Let $B=(-\alpha_1,\ldots,\alpha_n)$. If
$\lambda$ is an eigenvalue of $A^T\!A$ associated with an
eigenvector $\xi=(a_1,\ldots,a_n)^T\!$, i.e.\
$A^T\!A\xi=\lambda\xi$. Then $B^T\!B\eta=\lambda\eta$, where
$\eta=(-a_1,\ldots,a_n)^T\!$. Therefore $\lambda_{min}$, as
$\delta,\varphi$ and $\psi$, is independent of the choose of the
polyhedral cone $Q$ and is indeed determined by the hyperplanes
$H_1,\ldots,H_n$.
\end{rem}

\section*{Acknowledgements}
The author would like to thank professor Gu,~C.H. and Zhou,~Z.X. for
many questions, comments and suggestions. Thanks also go to
professor Hu,~H.S. for her staunch support.

\end{document}